\documentclass{amsart}
\usepackage{amssymb}
\usepackage{amsfonts}
\usepackage{latexsym}

\newtheorem{theorem}{Theorem}[section]
\newtheorem{lemma}[theorem]{Lemma}
\newtheorem{proposition}[theorem]{Proposition}
\newtheorem{corollary}[theorem]{Corollary}
\theoremstyle{definition}
\newtheorem{definition}[theorem]{Definition}

\newtheorem{remark}[theorem]{Remark}

\newcommand{\Hom}{{\rm {Hom}}}
\newcommand{\Ext}{{\rm {Ext}}}

\newcommand{\Aut}{\rm {Aut}}

\newcommand{\Rep}{{\rm {Rep}}}

\newcommand{\raro}{\to}

\newcommand{\ot}{\otimes}

\newcommand{\ben}{\begin{enumerate}}
\newcommand{\een}{\end{enumerate}}

\newcommand{\Vect}{{\rm {Vect}}}
\newcommand{\Lie}{{\rm {Lie}}}
\newcommand{\Cl}{{\rm {Cl}}}

\begin{document}

\title{On families of triangular Hopf algebras}

\author{Pavel Etingof}
\address{Department of Mathematics, Massachusetts Institute of Technology,
Cambridge, MA 02139, USA}
\email{etingof@math.mit.edu}

\author{Shlomo Gelaki}
\address{Department of Mathematics, Technion-Israel Institute of
Technology, Haifa 32000, Israel}
\email{gelaki@math.technion.ac.il}
\maketitle

%%%%%%%%%%%%%%%%%%%%%%%%%%%%%%%%%%%%%%%%%%%%%%%%%%%%%%%%%%%%%%%%%%%%%%%%%%%%%%%%%%%%%%%%%%%%%%%%%%%%%%%%%%%%%%%%%%%%%%%%%%%%%%%%%%%%%%%%%%%%%%%%%%%%%%%%%%%%

\begin{section}
{Introduction}
It was recently discovered ([AEG], [EG1]) that
the twisting procedure of finite supergroups can be used
to construct new families
of triangular Hopf algebras. Furthermore, it
turned out that these families
exhaust all triangular Hopf algebras
with the Chevalley property (in particular, all previously
known triangular Hopf algebras).
The goal of this paper is to continue the study of these families, and
in particular to give some applications outside of the theory of triangular
Hopf algebras.

More specifically, we consider the family of finite-dimensional complex
triangular Hopf algebras
$A(G,V,u,B)$, where $G$ is a finite group, $V$ a finite-dimensional
representation of $G$, $u$ a central element of order $2$ in $G$ acting
by $-1$ on $V$, and $B$ an element in $S^2V$. This is the simplest of the
families from [AEG], [EG1], corresponding to the twists which
are entirely contained in the ``superpart'' of a finite supergroup.

We start by finding the condition under which two members of such a family
are isomorphic as Hopf algebras (without regard for the triangular structure).
This allows us to calculate the moduli space of isomorphism
classes, which often turns out to be of positive dimension.
In particular, we see that this construction easily
produces continuous families of pairwise non-isomorphic Hopf algebras,
and hence Kaplansky's 10th conjecture [K] that there are finitely many
isomorphism classes of Hopf algebras in each dimension
(which was disproved in [AS], [BDG], [G]) fails even for triangular Hopf
algebras.
In fact, the lowest dimension in which we get continuous triangular families
is $32,$ which is the lowest dimension in which continuous families of Hopf
algebras are known to exist [Gr]. More precisely, in dimension $32$ we get
three non-equivalent $1-$parameter families, which are exactly the duals to
the three families constructed in [Gr] (and shown to be the only
continuous families of $32-$dimensional pointed Hopf algebras [Gr]). In particular,
we see that the families in [Gr] are cotriangular.

We also consider the question when the Hopf algebra $A(G,V,u,B)$ is twist
equivalent to the Hopf algebra $A(G,V,u,B')$ by twisting of the multiplication.
We show that if $(S^2V)^G=0$ then
for each $B$, there are only finitely many Hopf algebra isomorphism
classes of $A(G,V,u,B')$ with this property. This implies that
whenever the set of isomorphism classes is
infinite (e.g. in the $32-$dimensional examples of [Gr]),
we obtain continuous families of Hopf algebras which are not
equivalent by a twist of multiplication. This disproves a weakened
form of Kaplansky 10th conjecture suggested by Masuoka [M],
which states that such a family cannot exist.

Next, we study the algebra structure of $A(G,V,u,B)^*$,
and show that it is isomorphic as an algebra to
a direct sum of Clifford algebras. The only invariants
of Clifford algebras being the dimension and the rank of the quadratic form,
we see that in each family $A(G,V,u,B)$ there are finitely many
coalgebra types. Since all the members of this family
also have the same algebra type, one might wonder if
there are finitely many algebra (or, equivalently, coalgebra)
types of Hopf algebras
of a given dimension. We do not conjecture this, but it
would be interesting to find a counterexample.

Finally, we apply the results of this paper to the theory
of tensor categories. Namely, we start by recalling
Schauenburg's theorem that if $H_1,H_2$ are finite-dimensional
Hopf algebras such that $\Rep(H_1)$ is equivalent to $\Rep(H_2)$
as a  tensor category then $H_1,H_2$ are twist equivalent
(we give a simple proof of this result, which also works in the more general
case of quasi-bialgebras). This theorem implies that
if $(S^2V)^G=0$ then in the family of tensor categories $\Rep(A(G,V,u,B)^*)$,
each member is equivalent to at most finitely many other members.
On the other hand, there are only finitely many types of abelian
categories among them (since there are finitely many algebra types
of the underlying Hopf algebras). Thus we get an example of a {\it finite} abelian
category (i.e. one equivalent to the representation category of a
finite-dimensional algebra) which admits infinitely many
non-equivalent rigid
tensor structures with a fixed Grothendieck ring. 
As far as we know, such examples were not
previously available, and moreover it was shown by Ocneanu (unpublished)
that they cannot exist for semisimple categories.
\end{section}

\begin{section}
{The set of isomorphism classes of $A(G,V,u,B)$}
Throughout the paper, the ground field is the field of complex numbers
$\mathbb{C}.$

Let $G$ be a finite group, $V$ a finite-dimensional representation of
$G$ over $\mathbb{C},$ $u\in G$ a central element satisfying $u^2=1$ and
$u_{|V}=-1,$ and $B\in S^2V$
a symmetric $2-$tensor. Recall from [EG1] that a triangular Hopf algebra
$A=A(G,V,u,B)$ with the Chevalley property can be associated to $(G,V,u,B)$ in the
following way.

Let $\mathbb{C}[G\ltimes V]$ be the group algebra of the supergroup $G\ltimes
V$ (i.e., the Hopf superalgebra $\mathbb{C}[G]\ltimes \Lambda V$), and set
$J:=e^B\in \mathbb{C}[G\ltimes V]^{\ot 2}.$ Then it was shown in [AEG] that $J$ is
a twist for $\mathbb{C}[G\ltimes V],$ so we can form the Hopf superalgebra
$\mathbb{C}[G\ltimes V]^J,$
and then modify it using $u$ to get the Hopf algebra $A=A(G,V,u,B)$ (see [AEG,
Section 3]).

\begin{proposition}\label{isomthas}
Let $A_1=A(G_1,V_1,u_1,B_1)$ and $A_2=A(G_2,V_2,u_2,B_2),$ with
$V_1,V_2\ne 0.$ Then the Hopf algebras $A_1,A_2$ are isomorphic if and only if
there exists a group isomorphism
$\phi :G_1\to G_2$ such that $\phi(u_1)=u_2,$ and a $G_1-$isomorphism $\eta:V_1\to
\phi^*(V_2)$ such that $(\eta\ot \eta)(B_1)=B_2$ modulo $(S^2V_2)^{G_2}.$
\end{proposition}

\begin{proof}
Let us first prove the ``if'' direction. We need to show that $A(G,V,u,B)\cong
A(G,V,u,B+B')$ if $B'\in (S^2V)^{G}.$ To this end we note that for any $g\in
G,\alpha\in \Lambda V,$
and $J=e^B$ the coproduct in $\mathbb{C}[G\ltimes V]^J$ is given by
\begin{equation}\label{eq1}
\Delta^J(g\alpha)=
J^{-1}\Delta(g)\Delta(\alpha)J=
\Delta(g)\Delta(\alpha)e^{B-B^g},
\end{equation}
where $B^g(x,y):=B(gx,gy)$ (here we used the fact that $B$ is even, hence central
in $\Lambda V\ot \Lambda V$).
Therefore if we replace $B$ by $B+B',$ where $B'\in (S^2V)^{G},$ the
isomorphism type of the Hopf algebra $A(G,V,u,B)$ does not change.

Conversely, suppose $\xi:A_1\to A_2$ is an isomorphism of Hopf algebras. Taking
quotients by the radicals, we get an isomorphism of Hopf algebras
$\mathbb{C}[G_1]\to
\mathbb{C}[G_2],$ hence an isomorphism of groups $\phi:G_1\to G_2.$ Since the only
non-trivial $1:g$ skew-primitive elements in $A_1,A_2,$ arise for $g=u_1,u_2,$
respectively, we must have $\phi(u_1)=u_2.$ Therefore, $\xi$ induces an isomorphism
between the Hopf superalgebras $\mathbb{C}[G_1\ltimes V_1]^{J_1}$ and
$\mathbb{C}[G_2\ltimes V_2]^{J_2},$ hence between the spaces $V_1,V_2$ of their
primitive elements. Let us call this isomorphism $\eta.$ Clearly, $\eta$ is a
$G_1-$isomorphism $V_1\to \phi^*(V_2).$ Finally, $\xi$ induces an isomorphism
$\mathbb{C}[G_1\ltimes V_1]\to \mathbb{C}[G_2\ltimes V_2]^{J_2(\eta\ot
\eta)(J_1)^{-1}}.$ But this implies that $\mathbb{C}[G_2\ltimes V_2]^{J_2(\eta\ot
\eta)(J_1)^{-1}}$ is cocommutative, so
$J_2(\eta\ot \eta)(J_1)^{-1}=e^B,$ where $B$ belongs to $(S^2V_2)^{G_2}.$
\end{proof}

Let us denote by $\Aut(G,V)$ the group of all pairs $(\phi,\eta)$ where $\phi:G\to
G$ is an automorphism and $\eta:V\to \phi^*(V)$ is a $G-$isomorphism.

\begin{corollary}\label{c1}
The set of Hopf algebra isomorphism classes of $A(G,V,u,B),$ for fixed $G,V,u,$ is
$\displaystyle{\frac{S^2V/(S^2V)^G}{\Aut(G,V)}}.$
\end{corollary}

\begin{remark}
It follows from Corollary \ref{c1} that a non-trivial
continuous family of Hopf
algebras $A(G,V,u,B),$ with fixed $G,V,u,$ exists if and only if the quotient space
$S^2V/(S^2V)^G$ modulo the action of the group $\Aut(G,V)$ is infinite;
for example when $\dim(S^2V/(S^2V)^G)>\dim(\Aut(G,V)).$
\end{remark}

\end{section}

\begin{section}
{twist equivalence classes of $A(G,V,u,B)$}
Recall that for any $B,B'$, the Hopf algebras $A(G,V,u,B)$ and
$A(G,V,u,B')$ are twist equivalent as
triangular Hopf algebras (see [AEG]).
To the contrary, for the dual Hopf algebras we have the following theorem.

\begin{theorem}\label{main}
Let $G,V,u,B$ be as before, and suppose $(S^2V)^G=0.$ Then there exist only
finitely many $B'\in S^2V/\Aut(G,V)$ such that $A(G,V,u,B)^*$ is twist equivalent
to $A(G,V,u,B')^*.$
\end{theorem}

The rest of the section is devoted to the proof of Theorem \ref{main}.

Let $F$ be a functor and let $R^iF$ denote its $i-$th derived functor
($R^0F:=F$). Recall that by a result of Grothendieck, if $Q$ is an exact
functor then $R^i(QF)=QR^iF.$ Recall also that $R^i\Hom(*,N)=\Ext^i(*,N).$

In the following lemma we calculate the Hochschild
cohomology of $\mathbb{C}[G\ltimes V]$ with coefficients in the trivial bimodule
$\mathbb{C}$ (i.e.
$abc:=\varepsilon(a)\varepsilon(c)b,$ $a,c\in \mathbb{C}[G\ltimes V],$
$b\in \mathbb{C}$).
\begin{lemma}\label{cohom}
$H^i(\mathbb{C}[G\ltimes V],\mathbb{C})=(S^iV^*)^G.$
\end{lemma}
\begin{proof}
Let $F:\Rep(\mathbb{C}[G\ltimes V])\to \Rep(\mathbb{C}[G])$ be the
functor defined by $F(X)=\Hom_{\Lambda V}(X,\mathbb{C}),$ and let
$Q:\Rep(\mathbb{C}[G])\to \Vect$ be the functor defined by $F(Y)=Y^G.$
Then $Q$ is exact, and $QF(X)=\Hom_{\mathbb{C}[G]\ltimes \Lambda V}(X,
\mathbb{C}).$ Therefore, on the one hand,
$$R^i(QF)(\mathbb{C})=\Ext^i_{\mathbb{C}[G]\ltimes \Lambda
V}(\mathbb{C},\mathbb{C})=H^i(\mathbb{C}[G]\ltimes \Lambda V,\mathbb{C}),$$
and on the other hand,
$$QR^iF(\mathbb{C})=\Ext^i_{\Lambda V}(\mathbb{C},\mathbb{C})^G=
(H^i(\Lambda V,\mathbb{C}))^G=(S^iV^*)^G,$$
where the last equality follows from the Koszul duality.
\end{proof}

Recall that a linear form $\Phi:H\ot H\raro \mathbb{C}$
is called a {\em Hopf $2-$cocycle} for $H$ if it has an inverse
$\Phi^{-1}$ under the convolution product $*$ in $\Hom_{\mathbb{C}}(H\ot
H,\mathbb{C})$, and satisfies:
$$\sum \Phi(a_1b_1,c)\Phi (a_2,b_2)= \sum \Phi
(a,b_1c_1)\Phi (b_2,c_2)$$
for all $a,b,c\in H.$

Given a Hopf $2-$cocycle $\Phi$ for $H,$ one can construct a new Hopf
algebra
$$(H_{\Phi}, m_{\Phi},1, \Delta, \varepsilon,S_{\Phi})$$
as follows. As a coalgebra, $H_{\Phi}=H.$ The new multiplication is given by
$$m_{\Phi}(a\ot b)=\sum \Phi^{-1} (a_1,b_1)a_2b_2\Phi
(a_3,b_3)$$
for all $a,b\in H.$ The new antipode is given by
$$S_{\Phi}(a)=\sum \Phi ^{-1}(a_1,S(a_2))S(a_3)\Phi
(S(a_4),a_5)$$
for all $a\in H.$

Note that it is straightforward to verify that if $H$ is finite-dimensional,
then $\Phi\in H^*\ot H^*$ is a Hopf $2-$cocycle for $H$ if and only
if it is a twist for $H^*$ (here we do not impose the counit property on the
twist). So, in particular, the group $(H^*)^{\times}$ of
invertible elements in $H^*$ acts as gauge transformations
on the set of Hopf $2-$cocycles for $H.$

\begin{proposition}\label{h2c}
Let $H$ be any Hopf algebra which is isomorphic to $\mathbb{C}[G]\ltimes \Lambda V$
as an algebra, with the usual counit. Then there exist only finitely many
Hopf $2-$cocycles $\Phi$ on $H,$ modulo gauge transformations,
such that $H_{\Phi}$ is isomorphic to $\mathbb{C}[G]\ltimes \Lambda
V$ as an algebra, with the usual counit.
\end{proposition}
\begin{proof}
Let $X$ be the space of all Hopf $2-$cocycles such
that $H_{\Phi}\cong \mathbb{C}[G]\ltimes \Lambda V$ as algebras, with the usual
counit. Let $Y$ be the space of all Hopf $2-$cocycles on $H.$
Then $Y$ is an affine algebraic variety and $X\subseteq Y.$
Let $L:=(H^*)^{\times}$ be the group of invertible elements in $H^*;$ it acts on
$Y$ by gauge transformations and preserves $X.$ So we have a map
$\tau_x:\Lie(L)\to
T_xY,$
for all $x\in Y.$ It is easy to see from the definition of Hopf $2-$cocycles
that $T_xY=Z^2(H_x,\mathbb{C})$ (Hochschild
$2-$cocycles) and
$\tau_x(\Lie(L))=B^2(H_x,\mathbb{C})$ (Hochschild $2-$coboundaries).
Therefore,
$$T_xY/\tau_x(\Lie(L))=H^2(H_x,\mathbb{C})=H^2(\mathbb{C}[G]\ltimes
\Lambda V,\mathbb{C})=(S^2V^*)^G=((S^2V)^G)^*=0$$
using Lemma \ref{cohom}.
Thus, $\tau_x$ is surjective, which implies that $X/L$ is finite by Proposition
\ref{pap}.
\end{proof}

We are now ready to prove Theorem \ref{main}.

\noindent
{\it Proof of Theorem \ref{main}.}
By Proposition \ref{h2c}, there exist finitely many gauge equivalence classes of
Hopf $2-$cocycles $\Phi$ such that $A(G,V,u,B)_{\Phi}$ is isomorphic to
$\mathbb{C}[G\ltimes V]$ as an algebra, with the usual counit. This
implies that there exist
finitely many isomorphism classes of Hopf algebras $A(G,V,u,B')$ which are twist
equivalent to $A(G,V,u,B)$ by twisting of multiplication. By Corollary \ref{c1},
this implies the theorem.
\end{section}

\begin{section}
{$32-$dimensional examples}
As we remarked in the introduction,
examples of non-trivial continuous
families of triangular Hopf algebras
occur already in dimension $32.$
Here we describe these families, and identify them with
the families studied in [Gr].

More specifically,
in [Gr] Gra\~{n}a classified all complex $32-$dimensional pointed Hopf algebras
$H,$ and showed that they fall into finitely many isomorphism classes, except
for three $1-$parameter continuous families
(a family with $G(H)=\mathbb{Z}_4\times \mathbb{Z}_2$, and
two families with $G(H)=\mathbb{Z}_8$). We will prove the following result.

\begin{theorem}\label{examples}
The three continuous families of $32-$dimensional pointed Hopf algebras
of [Gr] are of the form $A(G,V,u,B)^*$ for appropriate
$G,V,u,B$.
\end{theorem}
\begin{proof}
To prove the theorem,
we will construct three families of $32-$dimensional pointed Hopf algebras
as $A(G,V,u,B)^*$, after which it will not be difficult to see
that they are identical to the families of [Gr].

\noindent
{\bf 1.} Let $G:=\mathbb{Z}_8=<a>,$ and
$\chi:\mathbb{Z}_8\to \mathbb{C}^*$ be the
character defined by $\chi(a^m)=e^{2\pi im/8}.$ Let $V:=\chi\oplus \chi^3=
sp\{e_1\}\oplus sp\{e_2\},$ where $a\cdot {e_1}=\chi(a)e_1$ and $a\cdot
{e_2}=\chi(a^3)e_2,$ and let
$u:=a^4$ (so $u_{|V}=-1$). Then any $B\in S^2V$ can be represented by a symmetric
$2\times 2$ matrix with respect to the basis $(e_1,e_2).$ We claim that
$A(G,V,u,B_1),A(G,V,u,B_2)$ are isomorphic if and only if $B_1=B_2$ in
$S^2V/\mathbb{Z}_2\ltimes D,$ where $\mathbb{Z}_2$ is generated by the matrix
$\left(\begin{array}{cc}
0 & 1\\
1 & 0
\end{array} \right),$
$D:=\Aut_G(V)$ is the group of invertible $2\times 2$ diagonal matrices, and the
action is given by $d\cdot B:=dBd,$ for all $d\in \mathbb{Z}_2\ltimes D.$
So, the set of isomorphism classes is $S^2V/\mathbb{Z}_2\ltimes D.$ Indeed,
$S^2V=\chi^2\oplus \chi^4\oplus
\chi^6,$ so $(S^2V)^G=0.$ Also the
isomorphism $\phi:G\to G$
is either determined by $a\mapsto a$ or by $a\mapsto a^3,$ so $\eta$ is
represented
either by a matrix of the form
$\left(\begin{array}{cc}
a_1 & 0\\
0 & a_2
\end{array} \right)$
for $a_1,a_2\ne 0,$ or by a matrix of the form
$\left(\begin{array}{cc}
0 & b_1\\
b_2 & 0
\end{array} \right)$
for $b_1,b_2\ne 0.$ Thus, the claim follows from Corollary \ref{c1}.

\noindent
{\bf 2.} Let $G,\chi,u$ be as in {\bf 1}, and let $V:=\chi\oplus
\chi^5 = sp\{e_1\}\oplus sp\{e_2\},$ where $a\cdot {e_1}=\chi(a)e_1$ and
$a\cdot {e_2}=\chi(a^5)e_2$ (so $u_{|V}=-1$).
As before, we claim that the set of isomorphism classes is
$S^2V/\mathbb{Z}_2\ltimes D.$
Indeed, in this case $S^2V=2\chi^2\oplus \chi^6,$ so as before
we have $(S^2V)^G=0,$ and the
result follows in a similar way to case {\bf 1}.

\noindent
{\bf 3.} Let $G=\mathbb{Z}_4\times \mathbb{Z}_2=<a>\times <b>.$ Let
$\chi:\mathbb{Z}_4\to \mathbb{C}^*$ be the character defined by
$\chi(a^m)=e^{2\pi im/4},$ $\chi_+:\mathbb{Z}_2\to \mathbb{C}^*$ the trivial
representation, and $\chi_-:\mathbb{Z}_2\to \mathbb{C}^*$ the non-trivial
representation. Let $V:=(\chi,\chi_+)\oplus (\chi,\chi_-),$ and $u:=(a^2,1).$
Then $u_{|V}=-1.$ We claim that, as in the previous two cases,
the set of isomorphism classes
is $S^2V/\mathbb{Z}_2\ltimes D$.
Indeed, we have $S^2V=2(\chi^2,\chi_+)\oplus
(\chi^2,\chi_-),$ so again $(S^2V)^G=0.$ Also,
it is straightforward to check that in this case the only
allowed
non-trivial automorphism $\phi:G\to G$ is the one determined by $a\mapsto
a,$ $b\mapsto
ba^2,$ which interchanges $(\chi,\chi_+)$ with $(\chi,\chi_-).$
Hence the result follows in a similar way to cases {\bf 1,2}.

Now, comparing this with [Gr], it is easy to see
that the above three families of Hopf algebras $A(G,V,u,B)^*$ exactly coincide with
the families of [Gr].
In detail, the family in example {\bf 3}
corresponds to the family in [Gr] with
$G(H)=\mathbb{Z}_4\times \mathbb{Z}_2$, and the families of examples
{\bf 1} and {\bf 2} correspond to the two families in [Gr] for
$G(H)=\mathbb{Z}_8$ (Line 2 and Line 4 of Table 15, respectively).
More precisely, it is not hard to check
that under this correspondence,
the lifting parameters $(\lambda_1,\lambda_2,\lambda_3)$
of [Gr] are related to our parameter $B$ via the formula
$B=\left(
\begin{array}{cc}
\lambda_1 & \lambda_3/2\\
\lambda_3/2 & \lambda_2
\end{array}
\right).$
The theorem is proved.
\end{proof}
\begin{remark} The equivalence relation between the lifting
parameters $(\lambda_1,\lambda_2,\lambda_3)$ is thus defined
by $(\lambda_1,\lambda_2,\lambda_3)\sim (t^2\lambda_1,s^2\lambda_2,
ts\lambda_3)$, $t,s\ne 0$, and
$(\lambda_1,\lambda_2,\lambda_3)\sim (\lambda_2,\lambda_1,\lambda_3)$.
This is exactly the equivalence relation which occurs in [Gr] (see the web version
of [Gr]; math.QA/0110033).
\end{remark}

\begin{corollary}
The $32-$dimensional examples $A(G,V,u,B)^*$
fall into infinitely many twist equivalence classes.
\end{corollary}
\begin{proof}
Follows from Theorems \ref{main},\ref{examples}.
\end{proof}
\end{section}

\begin{section}
{the algebra structure of $A(G,V,u,B)^*$}
In this section we describe the  algebra structure of the Hopf algebra
$A(G,V,u,B)^*.$ We start by recalling the definition of a Clifford algebra.

\begin{definition}
Let $V$ be a finite-dimensional vector space, and $B$ a symmetric bilinear form on
$V.$ The Clifford algebra $\Cl(V,B)$ is the quotient
of the tensor algebra $T(V)$ by the ideal generated by all elements of the form
$vw+wv-2B(v,w)1,$ for all $v,w\in V$ (e.g. $\Cl(V,0)=\Lambda V$).
\end{definition}

Note that the algebra $\Cl(V,B)$ has a unique structure of a superalgebra,
determined by requiring $V$ to be odd. Recall also that $\Cl(V,B)$ has a
filtration determined by letting $v\in V$ have degree $1,$ so that $gr
\Cl(V,B)=\Lambda V,$ and hence $\dim(\Cl(V,B))=2^{\dim(V)}.$

Finally, recall that if $V_1,V_2$ are finite-dimensional vector spaces, and
$B_1,B_2$ are symmetric bilinear forms on $V_1,V_2$ respectively, then
$\Cl(V_1\oplus V_2,B_1\oplus B_2)=\Cl(V_1,B_1)\ot \Cl(V_2,B_2),$
where $\ot $ denotes the tensor product of superalgebras.

The main result of this section is the following theorem.
\begin{theorem}\label{algstr}
$A(G,V,u,B)^*\cong \bigoplus _{h\in G/<u>} \Cl\left(V^*\bigoplus \mathbb{C},
\left(\begin{array}{cc}
B-B^h & 0\\
0 & 1
\end{array} \right)
\right),$
\linebreak where $h=\{g,gu\},$ and $B^h$ denotes $B^g=B^{gu}.$
\end{theorem}
In particular, Theorem \ref{algstr} implies the following.
\begin{corollary}
The Hopf algebra $A(G,V,u,B)$ is pointed if and only if $B\in (S^2V)^G;$ i.e.,
$A(G,V,u,B)\cong A(G,V,u,0)$ as Hopf algebras.
\end{corollary}

\begin{corollary}\label{fat}
The algebras $A(G,V,u,B)^*,$ with fixed $G,V,u,$ fall into finitely many
isomorphism classes.
\end{corollary}

The rest of the section is devoted to the proof of Theorem \ref{algstr}.

\begin{proposition}\label{p1}
Let $V$ be a finite-dimensional vector space, and let $B\in S^2V.$
Define the supercoalgebra $A_{V,B}$ to be $\Lambda
V$ as a vector space, with comultiplication $\Delta_B$ given by
$\Delta_B(\varphi)=\Delta(\varphi)e^B,$ $\varphi\in \Lambda V,$ where $\Delta$
denotes the usual comultiplication in $\Lambda V.$ Then, $A_{V,B}^*\cong
\Cl(V^*,B)$ as superalgebras.
\end{proposition}

\begin{remark}
By a super(co)algebra we mean a (co)algebra with an action of $\mathbb{Z}_2,$ and
the dual algebra is taken in the usual (rather than super) sense.
\end{remark}

\begin{proof}
It is clear that $A_{V_1\oplus V_2,B_1\oplus B_2}=A_{V_1,B_1}\ot
A_{V_2,B_2},$ where $\ot$ denotes the tensor product of supercoalgebras. Since
Clifford algebras have a similar property, and since any symmetric bilinear form
is a direct sum of $1-$dimensional forms, it is sufficient to prove the result when
$\dim(V)=1.$ But then the result is easy.
\end{proof}

Let us now determine the algebra structure of the superalgebra $A^*,$ where
\linebreak $A:=\mathbb{C}[G\ltimes V]^J$ and $J:=e^B.$
First, note that $A=\bigoplus_{g\in G} A_g,$ where
$A_g:=g\Lambda V,$ so $A^*=\bigoplus_{g\in G} A_g^*.$

\begin{proposition}\label{p2}
The following hold:
\ben
\item $A_g$ is a subcoalgebra for all $g\in G,$ so $A_g^*$ is a subalgebra for all
$g\in G$ and $A^*=\bigoplus_{g\in G} A_g^*$ as an algebra.
\item $A_g^*\cong \Cl(V^*,B-B^g),$ as superalgebras.
\een
\end{proposition}
\begin{proof}
The first part is clear from equation (\ref{eq1}), and
the second part follows from Proposition \ref{p1}.
\end{proof}

Assume now that we have fixed a central element $u\in G$ such that $u^2=1$ and
$u_{|V}=-1.$ Then we can consider the Hopf algebra $A(G,V,u,B).$ We now wish to
study the algebra structure of $A(G,V,u,B)^*.$

Recall that if $(C,\Delta)$ is a supercoalgebra, then one can define a coalgebra
$\overline{(C,\Delta)}=
(\overline C,\overline \Delta)$ such that $\overline C=\mathbb{C}[\mathbb{Z}_2]\ot
C=C\bigoplus uC,$ (i.e. $\mathbb{Z}_2$ is generated by $u,$ with $u^2=1$) and
$\overline
\Delta(ux)=(u\ot u)\Delta(x)$ and $\overline
\Delta(x)=\Delta_0(x)- (-1)^{p(x)}(u\ot 1)\Delta_1(x),$ for all $x\in C,$ where
$p(x)$ denotes the parity of $x,$ $\Delta_1(x)\in C\ot C_1$ and $\Delta_0(x)\in
C\ot C_0.$ Then, we have the following two straightforward results.

\begin{lemma}\label{l1}
$\overline C^*=\mathbb{C}[\mathbb{Z}_2]\ltimes C^*,$ where $\mathbb{Z}_2$
acts on $C^*$ by parity.
\end{lemma}

\begin{proposition}\label{p3}
Let $\widetilde {\Delta^J}$ be the coproduct of $A(G,V,u,B).$
Then, under $\widetilde{\Delta^J},$ $A_g\bigoplus A_{gu}$ is a subcoalgebra, for
all
$g\in G.$ Moreover, $(A_g\bigoplus A_{gu},\widetilde{\Delta^J})=\overline
{(A_g,\Delta^J)}.$
\end{proposition}
\begin{proof}
Follows from Proposition \ref{p2}(1).
\end{proof}

We can now prove Theorem \ref{algstr}.

\noindent
{\it Proof of Theorem \ref{algstr}.}
For any $h\in G/<u>,$ set $\widetilde{A_h}:=A_g\bigoplus A_{gu},$ where $g\in G$
represents $h.$ By Propositions \ref{p2}(2),\ref{p3} and Lemma \ref{l1},
$$\widetilde{A_h}^*=\mathbb{C}[\mathbb{Z}_2]\ltimes \Cl(V^*,B-B^h).$$
But this is by definition, equal to $\Cl\left(V^*\bigoplus \mathbb{C},
\left(\begin{array}{cc}
B-B^h & 0\\
0 & 1
\end{array} \right)
\right),$ as desired.

\end{section}

\begin{section}
{finite abelian categories with infinitely many tensor structures}
In this section we show that a finite abelian category (i.e. a category of the form
$\Rep(A),$ $A$ a finite-dimensional algebra)
may admit infinitely many non-equivalent tensor structures
with a fixed Grothendieck ring.

We start by formulating a theorem of Schauenburg. Although it will not be needed,
we will give a slightly generalized version of this theorem.
\begin{theorem}\label{tweq}
Two finite-dimensional quasi-bialgebras $H_1,H_2$ are twist
equivalent if and only if the categories $\Rep(H_1),\Rep(H_2)$ are tensor
equivalent.
\end{theorem}

\begin{proof}
Using the same proof as in Section 4.2 of [EG], one can show that if there exists a
tensor equivalence $\Rep(H_1)\to \Rep(H_2),$ which preserves ordinary dimensions of
objects, then $H_1,H_2$ are twist equivalent. So, it is sufficient to show that in
the finite-dimensional case, this property of a tensor equivalence is automatically
satisfied.

Consider the Grothendieck ring of $\Rep(H_1),$ with the distinguished basis formed by
the irreducible objects $V_1,\dots,V_n.$
Let $V\in \Rep(H_1),$ and consider the operator $L_V$ on the Grothendieck ring of
$\Rep(H_1),$ given by left multiplication by $V$ ($W\mapsto V\ot W$). Then $L_V$ is
represented by a matrix whose entries are non-negative integers. This
matrix has an eigenvector $(\dim(V_1),\dots,\dim(V_n)),$ with positive entries,
corresponding to the eigenvalue $\dim(V).$ But then, by the Frobenius-Perron
Theorem, $\dim(V)$ is the largest real eigenvalue of $L_V.$ Therefore, $\dim(V)$ is
preserved by any tensor equivalence, as desired.
\end{proof}

\begin{remark}
Theorem \ref{tweq} was proved in [S] (by a different method) under the assumption
that $H_1,H_2$ are
finite-dimensional Hopf algebras.
\end{remark}

As a consequence of Schauenburg's theorem we have the following.
\begin{theorem}
If $(S^2V)^G=0$ and $|S^2V/\Aut(G,V)|$ is infinite then for generic $B\in S^2V,$
the abelian category $\Rep(A(G,V,u,B)^*)$ has infinitely many distinct
rigid tensor structures, with the same Grothendieck ring as the tensor  
category $\Rep(A(G,V,u,B)^*)$.
\end{theorem}

\begin{remark}
For example, it is the case in the examples of Section 4 if
$\frac{\lambda_1\lambda_2}{\lambda_3^2}$ is a generic complex number.
\end{remark}

\begin{remark}
As we already mentioned in the introduction, Ocneanu showed that
such examples cannot exist for finite abelian {\it semisimple} categories. On the
other
hand, it is well known that there exist semisimple abelian categories
which are not finite (i.e. have infinitely many irreducible objects) and admit
infinitely many distinct tensor structures
with the same Grothendieck ring. For example, let $\mathfrak g$ be a
finite-dimensional complex semisimple Lie algebra. Then the representation
category of $U_q(\mathfrak g),$ for a generic complex number $q,$ is the same as
the representation category of $U(\mathfrak g),$ as an abelian category.
However, as tensor categories, these categories are distinct (up to the symmetry
$q\to q^{-1}$), although their Grothendieck rings are the same.
\end{remark}

\end{section}

\begin{section}
{appendix}
The argument used in the proof of Proposition \ref{h2c} relies on the following
standard proposition from algebraic geometry.

\begin{proposition}\label{pap}
Let $Y$ be an algebraic variety with an algebraic action of an
algebraic group $L.$ Let $Y_s$ be the subset of $Y$ consisting of the points $y$
where the natural linear map $\Lie(L)\to T_yY$ is surjective. Then $Y_s/L$ is finite.
\end{proposition}
Let us prove Proposition \ref{pap}, for the
convenience of the reader.

\begin{lemma}\label{help1}
Suppose $Z\subseteq Y$ are algebraic varieties
($Z$ is locally closed in $Y$), and
at a point $y\in Z,$ one has $T_yZ=T_yY.$ If $Z$ is smooth, there exists
a neighborhood $U$ of $y$ in $Y$ such that $U\subseteq Z.$
\end{lemma}
\begin{proof}
It is suffices to show that the inclusion $Z\to Y$
induces an isomorphism of local rings $g: O_y(Y)\to O_y(Z).$
Let $d:=\dim(T_yZ)=\dim(T_yY).$ It is clear that $\dim(O_y(Y))\ge d=\dim(T_yY),$ so
$y$ is a smooth point of $Y.$ The rest follows from the algebraic version of
the inverse function theorem.
\end{proof}

\begin{lemma}\label{help}
For any point $y\in Y_s,$ $Ly$ is open in $Y.$
\end{lemma}
\begin{proof}
Clearly, the lemma follows from Lemma \ref{help1}: we can take $Z:=Ly.$
\end{proof}

\noindent
{\it Proof of Proposition \ref{pap}.}
The proposition follows from Lemma \ref{help} since it implies that the orbits of
$L$ on $Y_s$ are disjoint open subsets of $Y,$ so the number of such orbits has to
be finite.
\end{section}

\begin{section}
{Acknowledgments}
We are grateful to S. Montgomery for
turning our attention to Masuoka's question, to M. Gra\~{n}a for
inspiring discussions and explanations of his work, and to V. Ostrik for
useful conversations.

P.E.\ was partially supported by the NSF grant DMS-9988796.
P.E.\ partially conducted his research for the
Clay Mathematics Institute as a Clay Mathematics Institute Prize Fellow.
S.G.\ thanks the Mathematics Department of MIT for the warm hospitality
during his visit. S.G.\ research was supported by the VPR - Fund at the
Technion and by the fund for the promotion of research at the Technion.
\end{section}

\bibliographystyle{ams-alpha}

\end{document}